\documentclass[11pt]{elsarticle}

\usepackage{lineno}
\usepackage[colorlinks=true,breaklinks=true,pdftex]{hyperref}
\usepackage{amsfonts,amsmath,amssymb,verbatim}
\usepackage{graphicx,color}
\usepackage{a4wide}
\usepackage{xfrac}
\usepackage{todonotes}

\biboptions{authoryear}

\newtheorem{dfn}{Definition}

\newtheorem{lem}[dfn]{Lemma}

\newtheorem{thm}[dfn]{Proposition}

\newtheorem{exm}[dfn]{Example}

\DeclareMathOperator{\spn}{Span}
\DeclareMathOperator{\im}{Im}

\graphicspath{{./Figures/}}

\modulolinenumbers[1]

\journal{CAGD}

\begin{document}

\begin{frontmatter}

\title{Lie Group Approach to Envelope Surfaces} 



\author{Michal Molnár$^a$}
\ead{mysymol@mmolnar.sk}
\address[1]{Charles University, Faculty of Mathematics and Physics\\Sokolovsk\'a 83, Prague 186 75, Czech Republic.}

\author{Zbyn\v{e}k \v{S}\'{\i}r$^{a}$}
\ead{zbynek.sir@mff.cuni.cz}

\author{Jana~Vr\'{a}bl\'{\i}kov\'a$^{b,*}$}
\ead{jana.vrablikova@inria.fr}
\address[2]{Centre Inria d’Universit\'e C\^{o}te d’Azur\\
2004 Route des Lucioles, B.P. 93, Sophia Antipolis, 06902, France.}
\cortext[cor]{Corresponding author}


%
%
\begin{abstract}
In this paper, we develop a new and efficient approach to the computation of envelope surfaces. We interpret one-parameter systems of surfaces as curves in the homogeneous spaces of suitable Lie groups. Using the formalism of Lie groups and Lie algebras, we rigorously capture the inherent symmetry and linearity in the computation of envelopes. In particular, the possible set of characteristic curves (which constitute the envelope surface) can be precomputed as the intersection of a fixed canonical surface and a low-dimensional set of its possible “derivatives.”

To demonstrate the effectiveness of our approach, we present several examples of surfaces undergoing transformations from various Lie groups. As a remarkable side result, we show that the characteristic curves and the envelopes of cones undergoing rational motions are themselves rational. Furthermore, we provide an explicit rational parameterization of these envelopes and use it to solve the trimming problem.
\end{abstract}


\end{frontmatter}



\section{Introduction} 

\label{intro}
The main topic of this paper is the connection between envelope surfaces and Lie groups. Lie groups and their homogeneous spaces are two of the most widely used mathematical concepts, see e.g. \cite{LieHomo, matrix,Karger} for an accessible introduction. They are universally applicable for expressing the symmetries that occur in virtually any scientific domain. Kinematics \cite{Karger} and mathematical physics \cite{phys}, in particular, are largely based on Lie group theory. To a smaller extent, they have also found applications in geometric modeling, see e.g. \cite{LieAppl,Wallner}. Envelope surfaces have been well studied for both their theoretical aspects and their applications. Indeed, they are considered a standard topic in differential geometry \cite{struikEnvelopes}. Envelopes are used as auxiliary constructions for various purposes, see e.g. \cite{kosinka,sir}. Various special envelope surfaces have been studied in geometric modeling, such as canal surfaces \cite{envdef} and developable surfaces \cite{dev}. The envelopes of moving (truncated) cones are of a great interest due to their applications in CNC machining \cite{barton,novyBarton}.

Our paper contributes to the research effort to find analytical expressions and parameterizations of the envelope surfaces, see e.g. \cite{Bi,flaquerenvelopes,peternell}. Although Lie groups are occasionally mentioned in this context, their structure is not systematically and generally exploited. In this regard, our results in Section \ref{Stheory} are original. The application to canal surfaces in Section \ref{apl} provides a new, compact proof of the well--known fact that they are rational as envelopes of a one--parameter system of spheres that depend rationally on the parameter \cite{envdef}. The rationality of the envelopes of moving cones is stated and proven in \cite{Bi,juettler1999,peternell,peternell2010,peternell2008, peternell2000}. In all cases, however, a dual approach is used, and no explicit formula is provided. In this paper, we provide an explicit rational parameterizations of these envelopes and use them to solve the trimming problem efficiently. 

The remainder of the paper is organized as follows: Section \ref{sec_Preliminaries} gathers the necessary definitions and properties of envelope surfaces, Lie groups, and their homogeneous spaces. In Section \ref{Stheory} we show that all surfaces of a certain fixed type naturally form a homogeneous space of a suitable Lie group. This space is then represented in a vector space of implicit functions. One--parameter systems of surfaces appear naturally as curves in the homogeneous space. We exploit the symmetry of the homogeneous space and the linearity of its tangent space to simplify the description of the characteristic curves. Section \ref{apl} applies this general approach to pipe and canal surfaces, freely changing ellipsoids, and, most importantly, cones undergoing Euclidean motions.  In the latter case, we provide rational parameterizations of the characteristic curves and the envelope and use it to solve the trimming problem. Finally, we conclude the paper.

\section{Preliminaries}\label{sec_Preliminaries}
This paper focuses on the interconnection between envelope surfaces and Lie groups. This section summarizes both topics.  In particular, we address envelope computation for one-parameter systems of implicitly described surfaces in $\mathbb R^3$. To better explain all the fundamental concepts, we provide an example that runs throughout the paper. Regarding Lie group theory, we first provide a general definition and then focus on the special Euclidean group $\mathbb{SE}(3)$, its Lie algebra, and its homogeneous spaces. 

\subsection{Envelope surfaces}

Let us consider a one--parameter system of surfaces $F_t$ embedded in $\mathbb R^3$ and indexed by a parameter $t$ over a real interval $I$, and their envelope surface $\chi$. While the previous sentence is intuitive, the key notions of {\em surface} and {\em envelope} are often not handled precisely in the research literature. By {\em surface}, a smooth $2$-dimensional manifold embedded in $\mathbb R^3$ is mostly meant. However, due to the usefulness of implicitly defined algebraic surfaces, sets of singular points (of measure zero) are typically permitted.

The intrinsic difficulty of the definition arises from two natural views on the {\em envelope}. The first view defines the envelope surface $\chi$ as the boundary of the union over $t$ of the volumes delimited by $F_t$ for each $t$. However, it is difficult to define these volumes in a general way and avoid unwanted overlaps. 
The second approach defines the envelope $\chi$ as the maximal set (surface) such that $\chi\cap F_t$ is one--dimensional for generic $t$, the two surfaces are tangent along this intersection curve, and a generic point of the envelope $\chi $ belongs to precisely one $F_t$. This approach is difficult to handle because of the non--generic points and the additional unwanted intersections of $F_t$ and the envelope.
In \cite{envdef}, a formally rigorous definition is given in the context of canal surfaces that can be generalized to algebraic implicit surfaces.

The strategy presented in this paper is relevant for the generally considered systems of surfaces and their envelopes. However, its computational aspect is most efficient for implicitly defined systems. For this reason, we will not elaborate on the definition of the envelope and, like many other authors, will rely on the following characterization, see e.g. \cite{flaquerenvelopes,struikEnvelopes}.

\begin{thm}\label{thm_mult}
Let a one--parameter system of surfaces be given by the implicit equations 
\begin{equation}\label{impl}
    F_t = \{[x,y,z]\in \mathbb R^3;~f(x,y,z,t) = 0\}, \qquad t\in I\ ,
\end{equation}

where $f(x,y,z,t)$ is a smooth function.  Then the characteristic curves are given as 
\begin{equation}\label{char}
    \chi_t=\{[x,y,z]\in \mathbb R^3; f(x,y,z,t) = 0 \land \frac{\partial f}{\partial t}(x,y,z,t) = 0 \}
\end{equation}

 and the envelope surface as their union $\displaystyle \chi=\bigcup_{t\in I}\chi_t$. 
\end{thm}

We now present an example of a system of quadratic surfaces that depends rationally on $t$. For an idea of the computational complexity involved, we provide not only figures, but also equations.

\begin{exm} \label{ex1}\upshape
    \begin{figure}[htb]
        \centering
        \includegraphics[height=5cm]{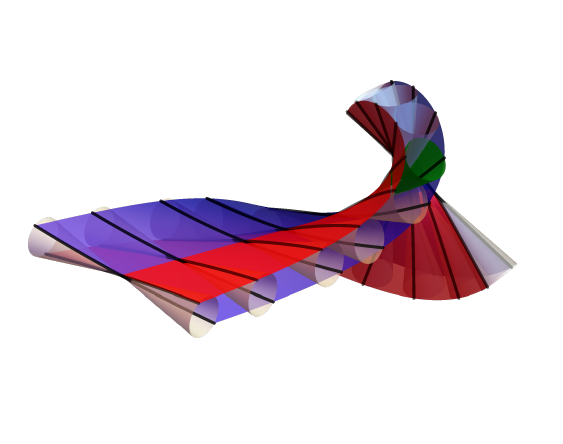}
        \caption{One parameter system of surfaces with their characteristics (black) and the envelope surface, which consists of two sheets (red and blue). The green cone shows the element of $F_t$ for $t = \frac 1 2$.
        }
        \label{fig:ex1}
    \end{figure}

Consider the one--parameter system of surfaces $F_t$, $t\in [0,1]$, given implicitly by \eqref{impl} with 
{\footnotesize
\begin{equation}\label{impF}
    \begin{split}
        f(x,&y,z,t)=\\ & x^2(107430976 t^4+233406592 t^3+157385824 t^2+10446368 t+4338116)\\
        &+y^2(108682304 t^4+195440768 t^3+186419296 t^2+68579872 t+21322564)\\
        &+z^2(24508864 t^4-86590592 t^3+2451616 t^2+80672992 t+26727964)\\
        &+xy(29392896 t^4-78061568 t^3-12087296 t^2+61155328 t+22004736)\\
        &+xz(77635584 t^4-45606912 t^3-291359744 t^2-125152768 t-289536)\\
        &+yz(-74400768 t^4-43902976 t^3+118962688 t^2+61441536 t+142272)\\
        &+x(-750809344 t^5-1016063680 t^4-311980800 t^3+121132640 t^2+43314320 t-23582412)\\
        &+y(-487677312 t^5-925983040 t^4-1463723200 t^3-1064846880 t^2-424846840 t-54543076)\\
        &+z(-79122048 t^5+589686720 t^4+561255360 t^3-47853600 t^2-245226120 t-79820724)\\
        &+1601008384 t^6+2264876032 t^5+2216961040 t^4+1726892960 t^3\\
        &+1173047560 t^2+482773672 t+91188121 \ ,
    \end{split}
\end{equation}
}
see Figure \ref{fig:ex1}. By the direct differentiating we get 
{\footnotesize 
    \begin{equation}
        \begin{split}
        \frac{\partial f}{\partial t}(x,&y,z,t)=\\ & x^2(429723904 t^3+700219776 t^2+314771648 t+10446368)\\
        &+y^2(434729216 t^3+586322304 t^2+372838592 t+68579872)\\
        &+z^2(98035456 t^3-259771776 t^2+4903232 t+80672992)\\
        &+xy(117571584 t^3-234184704 t^2-24174592 t+61155328)\\
        &+xz(310542336 t^3-136820736 t^2-582719488 t-125152768)\\
        &+yz(-297603072 t^3-131708928 t^2+237925376 t+61441536)\\
        &+8x(5414290 + 30283160 t - 116992800 t^2 - 508031840 t^3 - 469255840 t^4)\\
        &+8y(-304798320 t^4-462991520 t^3-548896200 t^2-266211720 t-53105855)\\
        &+8z(-49451280 t^4+294843360 t^3+210470760 t^2-11963400 t-30653265)\\
        &+8(1200756288 t^5+1415547520 t^4+1108480520 t^3+647584860 t^2+293261890 t+60346709) \ .
\end{split}
    \end{equation}
}
\begin{figure}[htb]
    \centering
    \includegraphics[height=5cm]{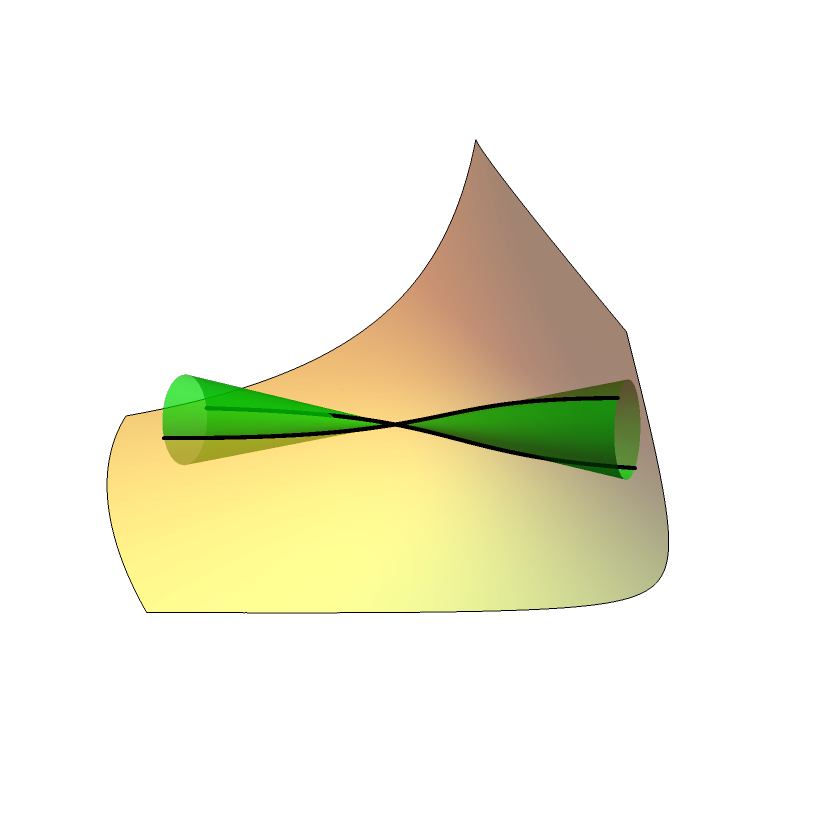}
    \caption{The surface $F_{t}$ for $t=\frac{1}{2}$ (green), the corresponding derivative surface (yellow), and the characteristic curve $\chi_t$ (black) as their intersection.}
    \label{fig:ex11}
\end{figure}

For a fixed $t\in[0,1]$, the characteristic curve $\chi_t$ on $F_t$ is obtained as the set of all points $[x,y,z]\in \mathbb R^3$ satisfying 
\begin{equation}
    f(x,y,z,t) =0 \quad \text{and} \quad \frac{\partial f}{\partial t}(x,y,z,t) = 0 \ ,\label{ex1f}
\end{equation}
see Figure \ref{fig:ex11}. This curve is a priori not rational. The envelope of the system (see Figure \ref{fig:ex1}) can most easily be obtained implicitly by eliminating the parameter $t$ from the equations \eqref{ex1f}. In this example, we were able to perform this elimination using the Gröbner basis technique. However, the resulting implicit equation is too long to include in the paper.
\end{exm}

Clearly, providing an efficient description of the envelope is not a trivial task. The implicit description is computationally involved and not entirely suitable, as it does not allow simple restriction to a given interval for $t$, nor efficient trimming of the surfaces $F_t$. For this reason, parameterizing the individual characteristic curves, at least approximately, is highly desirable. The following simple observation is useful for this purpose.
    
\begin{lem}\label{thm_envelopechar}
The characteristic curves given by the formula \eqref{char} are identical for the surfaces defined as the zero sets of the implicit functions $f(x,y,z,t)$ and the implicit functions $f^{\lambda}(x,y,z,t): = \lambda(t) f(x,y,z,t)$ where  $\lambda(t)$ is a smooth function that is non--zero for $t\in I$. Consequently, the envelope surfaces are identical as well. 
\end{lem}
    
\begin{prf}
The proof follows from direct computation. For any fixed values of $x,y,z\in \mathbb R$ and $t\in I$, the equations 
    \begin{eqnarray*}
        f^{\lambda}(x,y,z,t) &=& \lambda(t) f(x,y,z,t)=0 \ ,\\
        \frac{\partial f^{\lambda}}{\partial t}(x,y,z,t) &=& \lambda'(t)f(x,y,z,t) + \lambda(t)\frac{\partial f}{\partial t}(x,y,z,t)=0 
    \end{eqnarray*}
clearly define the same surfaces as the equations 
    \begin{equation*}
        f(x,y,z,t)=0, \quad
        \frac{\partial f}{\partial t}(x,y,z,t)=0\ ,
    \end{equation*}
because $\lambda(t)\neq 0$.\qed
\end{prf}

\medskip
While the previous Lemma merely confirms the intuition that any non--zero multiple of the respective equations equally well represents the implicitly given surfaces, it is interesting to observe that, by suitably setting the values of $\lambda(t)$ and $\lambda'(t)$  for a specific $t\in I$, the derivative surface $\frac{\partial f^\lambda}{\partial t}(x,y,z,t)=0$ (yellow in Figure \ref{fig:ex11}) can be any surface of the pencil given by the two surfaces $f(x,y,z,t)=0$ and $\frac{\partial f}{\partial t}(x,y,z,t)=0$. We plan to formally study such \emph{characteristic pencils} in our future research.

\subsection{Lie groups and their homogeneous spaces}

Let us now briefly review the essential concepts of Lie groups and their homogeneous spaces. This is a standard topic, and we refer the reader to \cite{LieHomo,Karger} for more details.

\begin{dfn}
A real smooth manifold $\mathcal G$ of dimension $n$ with a binary operation $\circ$ is called {\em Lie group} if
    \begin{enumerate}
        \item $(\mathcal G, \circ)$ is a group. 
        \item The mappings $\mathcal G\times \mathcal G \to \mathcal G: (x,y)\to x\circ y$ and $\mathcal G\to \mathcal G: x\to x^{-1}$ are smooth.          
    \end{enumerate}
The tangent space to $\mathcal G$ at its unit element $T_{\mathbf 1}\mathcal G$ is called the associated Lie algebra $\mathfrak g$. 
\end{dfn}

In other words, Lie groups combine the algebraic structure (that of a group) and the analytical structure (that of a smooth manifold) in a coherent way. The Lie algebra is equipped with a binary operation called the {\em Lie bracket}, which will not be used here. The most natural non--trivial examples of Lie groups are transformation groups, such as the group of all isometries, the group of all affine transformations, and the group of all conformal transformations, among others.

\begin{exm}\upshape
As a prime example, consider the group of all direct isometries of the Euclidean space $\mathbb R^3$, which is sometimes called the  {\em special Euclidean group}. This group has dimension 6 and is denoted $\mathbb{SE}(3)$. All the direct isometries are of the form 
    
\begin{equation}\label{eucl3}
\begin{pmatrix}x\\y\\z\end{pmatrix}\to A \begin{pmatrix}x\\y\\z\end{pmatrix}+\mathbf p \ ,
\end{equation}
where $A$ is a rotation $3\times 3$ matrix such that $AA^T=I$, $\det A=1$, and $\mathbf p \in \mathbb R^3$ is a translation vector. Thus, the pairs $( A, \mathbf p)$ represent all the elements of $\mathbb{SE}(3)$. The composition rule is given simply by composing the transformations:

\begin{equation*}
    \begin{pmatrix}x\\y\\z\end{pmatrix}\to B\left ( A \begin{pmatrix}x\\y\\z\end{pmatrix}+\mathbf p \right  ) +\mathbf q= B A \begin{pmatrix}x\\y\\z\end{pmatrix}+\left (B\mathbf p   +\mathbf q\right ) \ ,
\end{equation*}  
i.e., $(B, \mathbf q)\circ (A, \mathbf p)=(BA, B \mathbf p+\mathbf q)$. This group can be embedded into the group of regular $4\times 4$ matrices via 

\begin{equation*}
    \left( A, \mathbf p\right) \to 
    \left(\begin{array}{cccc}
    & & & \\
    &A&&\mathbf p\\
    & & & \\
    0&0&0&1
    \end{array} \right) \ .
\end{equation*}

The isometry in coordinates is then expressed as
\begin{equation*}
    \begin{pmatrix}x\\y\\z\\1 \end{pmatrix}\to \left (
    \begin{array}{cccc}
        & & & \\
        &A&&\mathbf p\\
        & & & \\
        0&0&0&1
    \end{array}
       \right )
    \begin{pmatrix}x\\y\\z\\1 \end{pmatrix} \ .
\end{equation*}
This embedding enables us to describe explicitly the corresponding Lie algebra
\begin{equation}\label{Lalg}
    \mathfrak g= \mathfrak{se}(3)=
    \left \{
    \begin{pmatrix}
        0&-a&-b&e\\
        a&0&-c&f\\
        b&c&0&g\\
        0&0&0&1
    \end{pmatrix}: 
    a,b,c,e,f,g \in \mathbb R \right \} \ ,
\end{equation}
see \cite{Karger} for more details. 
\end{exm}

As we will see later, the theory of homogeneous spaces provides the proper framework for describing systems of surfaces of a fixed type.

\begin{dfn}
    A non--empty set $\mathcal H$ is called a homogeneous space of a Lie group  $\mathcal G$ if $\mathcal G$ acts transitively on $\mathcal H$. More precisely, there is a group homomorphism $\Phi:\mathcal G\to Aut(\mathcal H)$ such that 
    \begin{equation}\label{tran}
        \forall x,y \in \mathcal H \quad \exists g \in \mathcal G : \Phi(g)(x)=y \ .
    \end{equation}      
\end{dfn}

In this context, $Aut(X)$ denotes either the group of all the bijections of $X$ or its subgroup consisting of bijections that preserve some structure of $X$. Thus, it contains only continuous, smooth, or other such bijections. The action mapping $\Phi$ is often omitted from the notation. That is, we write $g(x)$ instead of $\Phi(g)(x)$.

Unlike Lie groups, homogeneous spaces have no algebraic structure. Conversely, the action of the Lie group expresses the symmetry of the homogeneous space. There are no exceptional elements in $\mathcal H$ (such as the unit element); they are all equivalent because the group's action is transitive. However, it does have the structure of a smooth manifold, which it inherits from the Lie group. Due to the transitivity condition, we have
\begin{equation*}
    \mathcal H \simeq \mathcal G/\mathcal G_x \ ,   
\end{equation*}
where $\mathcal G_x$ is the subgroup of the stabilizers of some fixed $x\in \mathcal H$,
\begin{equation*}
    \mathcal G_x=\{g\in \mathcal G: g(x)=x\} \ .
\end{equation*}

As a trivial example, consider the affine space, which is the principal homogeneous space of the associated vector space. The term {\em principal} means that precisely one $g \in \mathcal G$ realizes the transitivity in \eqref{tran}.
    
The group $\mathbb{SE}(3)$ acts naturally on $\mathbb R^3$, making it its homogeneous space. The stabilizer of the point $(0,0,0)^T\in \mathbb R_3$ is the group of the orthonormal matrices with positive determinant, $\mathbb{SO}(3)$. Therefore, $\mathbb R_3 \simeq \mathbb{SE}(3)/\mathbb{SO}(3)$.

More importantly, the group $\mathbb{SE}(3)$ acts transitively on sets of mutually--isomorphic objects. For example, the set of all straight lines, the set of all planes, and the set of all spheres of a fixed radius are all homogeneous spaces of the group $\mathbb{SE}(3)$.    

\section{Systems of surfaces as curves in homogeneous spaces}\label{Stheory}
In this section, we show that all possible instances of a changing surface form a homogeneous space that can be embedded in a linear space of implicit equations. This embedding is a crucial step that allows us to fully exploit the symmetry and linearity in the computation of the characteristic curves.


\begin{dfn}
Let $\mathcal G$ be a Lie group acting on $\mathbb R^3$ and let $\mathcal H$ be its homogeneous space such that $\exists \bar F \in \mathcal H$ satisfying
\begin{equation*}
 \mathcal H= \{g(\bar{F}): g\in \mathcal G\} \ .   
\end{equation*}
Let $F_t$ be a one--parameter system of surfaces. We say that $F_t$ is \emph{contained} in $\mathcal H$ if $F_t\in \mathcal H$ for all $t \in I$. The surface $\bar F$ is called an \emph{elementary surface}. 
\end{dfn}

In other words, the homogeneous space $\mathcal H$ of a group $\mathcal G$ is the set of all transformations of one elementary surface $\bar F$ by elements of the group, and the group is large enough to produce the surfaces $F_t$ for all $t \in I$. Due to the nature of homogeneous spaces, any element of $\mathcal H$ can be taken as elementary. However, for computational purposes, we usually fix one particular surface $\bar F$. The choice of $\mathcal G$ and $\bar F$ so that $\mathcal H$ contains a particular one-parameter system of surfaces is usually straightforward, because we typically deal with some kind of transformation of a surface from the beginning.

\begin{exm}\label{ex2}\upshape
In Example \ref{ex1}, all surfaces in the system $F_t$ were constructed as cones of revolution with an angle $\alpha = \arctan(\frac{1}{5})$ between the axis and the straight lines (generators) of the surface. These surfaces are thus mutually related via direct isometries. In other words, system $F_t$ is contained in the homogeneous space $\mathcal H$ with the Lie group $\mathcal G=\mathbb{SE}(3)$ and a suitable elementary surface $\bar F$, say 
\begin{equation}\label{eq_elemCone}
    \bar F=\{[x,y,z]: x^2+y^2-\frac{1}{25}z^2 = 0\} \ .
\end{equation}
Let us stress that the surface $\bar F$ does not need (and in fact it does not) belong to the system $F_t$, see Figure \ref{fig:ex33}, left. 

\begin{figure}[htb]
    \centering
    \includegraphics[height = 5cm]{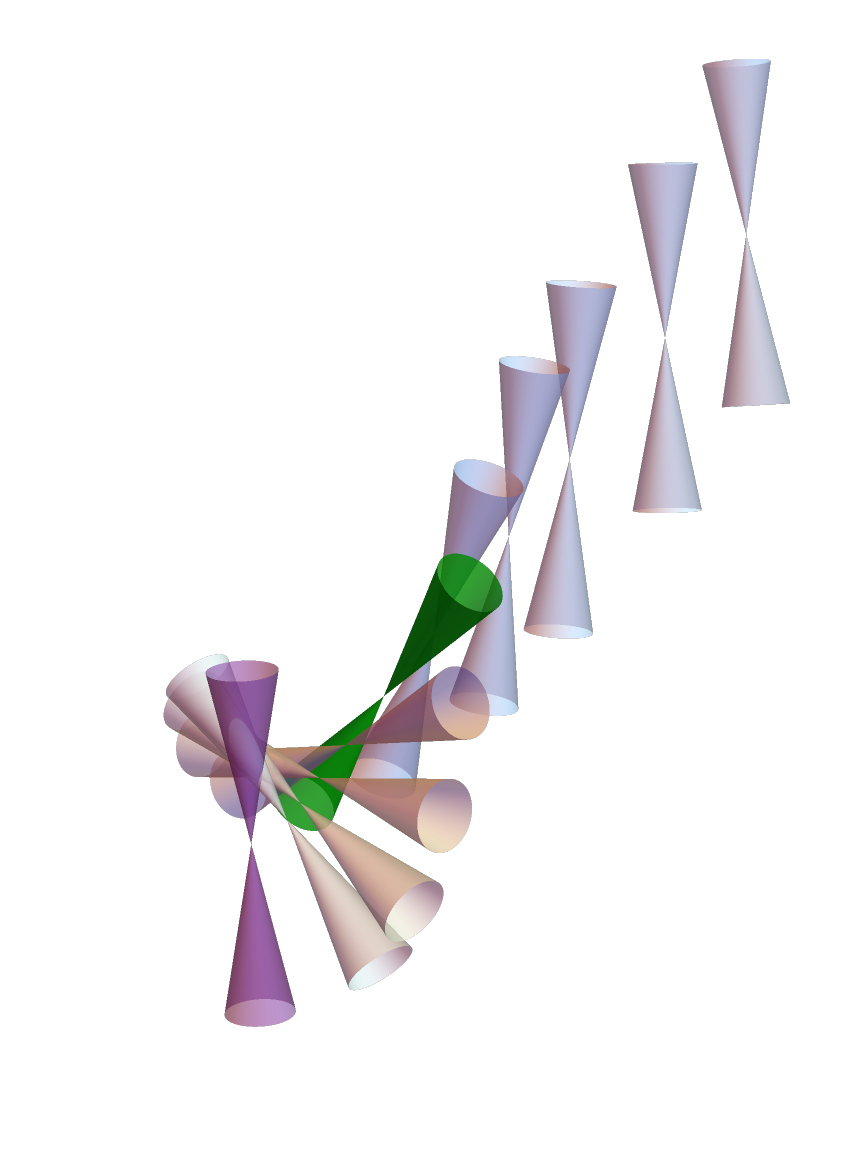}
\qquad \qquad \includegraphics[height = 5cm]{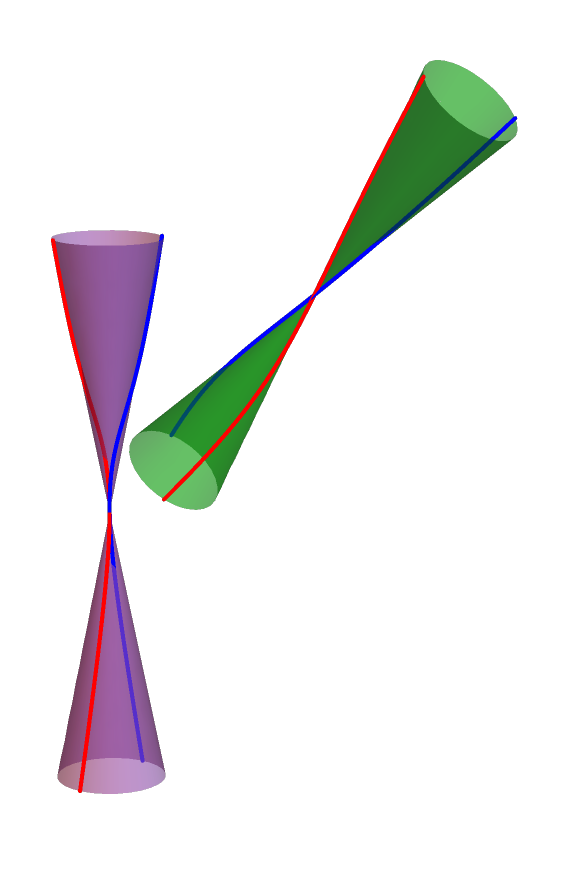}
 \caption{The mutual position of the one parameter system of cones and the elementary cone $\bar F$ (purple). The cone $F_t$ for $t=\frac{1}{2}$ is highlighted in green (left). Characteristic curves are transformations of curves on the elementary surface (right).}
    \label{fig:ex33}
\end{figure}
\end{exm}

In the previous example, we saw that the group of direct Euclidean isometries was large enough to contain the system $F_t$. However, if the system included general cones, a larger group (e.g., an affine group) would be necessary. For computational reasons, it is best to choose the smallest possible Lie group.

The next step is to describe a one--parameter system of surfaces solely within the group~$\mathcal G$.

\begin{dfn}
We say that a smooth parametric curve $c:I \to \mathcal G:t\to g_t$ is a {\em smooth lift} of the one--parameter system $F_t$ if $F_t=g_t(\bar{F})$ for all $t\in I$. 
\end{dfn}

A smooth lift is not unique because there are usually many elements $g\in \mathcal G$ such that $g(\bar F)=F_t$ for a $t \in I$. It does not need to exist globally, but the theory of homogeneous spaces shows that a lift exists locally. Moreover, this is usually not an issue because we often start with $g_t$ to describe the system $F_t$, e.g., in the context of CNC manufacturing. In fact, the system $F_t$ in Examples \ref{ex1} and \ref{ex2} was constructed by taking $F_t=g_t(\bar F)$, where 
\begin{equation*}
    g_t= 
\left(
\begin{array}{cccc}
 \frac{-468 t^2-788 t+123}{1108 t^2+788 t+517} & \frac{928 t^2+848 t+192}{1108 t^2+788 t+517} & \frac{384 t^2-736 t-464}{1108 t^2+788 t+517} &\frac 3 2 + 2t\\
 \frac{864 t^2+944 t+256}{1108 t^2+788 t+517} & \frac{588 t^2+268 t+387}{1108 t^2+788 t+517} & \frac{-368 t^2+272 t+228}{1108 t^2+788 t+517}&\frac 1 2 + 3t \\
 \frac{-512 t^2+608 t+432}{1108 t^2+788 t+517} & \frac{144 t^2-496 t-284}{1108 t^2+788 t+517} & \frac{-972 t^2-1292 t-3}{1108 t^2+788 t+517} & \frac 3 2 + 3 t\\
 0&0&0&1
\end{array}
\right)\in \mathbb{SE}(3)
\end{equation*}
is a smooth parametric curve in $\mathbb{SE}(3)$ and $\bar F$ is the cone \eqref{eq_elemCone}.

To fully exploit the differential structure of the Lie groups, the last ingredient required is an embedding of the homogeneous space $\mathcal H$ into a space of implicit functions.

\begin{dfn}
Let $\mathcal H$ be a homogeneous space of a a Lie group $\mathcal G$ such that $\mathcal H = \{g(\bar F): g \in \mathcal G\}$ for a $\bar F \in \mathcal H$. Let $\bar{f}(x,y,z)$ be an implicit equation of $\bar F$, i.e., $\bar F=\{[x,y,z]: \bar{f}(x,y,z) = 0\}$. We define the linear space $Q$ of the functions in three variables as 
\begin{equation*}
    Q=\spn\{f(x,y,z):=\bar{f}(g^{-1}(x,y,z)), g\in \mathcal G\}.
\end{equation*}
This space is called the {\em implicit representation space} of the homogeneous space $\mathcal H$. There is also a natural smooth mapping 
   \begin{align}     \label{repM}
   \phi:\ &\mathcal G \to Q \nonumber \\
   &g \mapsto \bar{f}\circ g^{-1}
\end{align}  
called the {\em implicit representation mapping}. 
\end{dfn}

The implicit representation space $Q$ is constructed to contain the implicit equations of all the surfaces in the homogeneous space $\mathcal H$. Indeed, $\mathcal H$ consists of surfaces $g(\bar{F})$ and clearly 
\begin{equation*}
    [x,y,z]\in g(\bar{F})  \Leftrightarrow g^{-1}(x,y,z) \in \bar{F}\Leftrightarrow\bar{f}(g^{-1}(x,y,z))=0 \ .
\end{equation*}
The implicit representation space may be infinite--dimensional for general surfaces. However, if $\bar F$ is algebraic, the space $Q$ has a finite dimension.

In Example \ref{ex2}, the elementary surface $\bar{F}$ is the zero set of the function $\bar{f}(x,y,z)=x^2+y^2-\frac{1}{25}z^2$. All its transformations $\bar{f}(g^{-1}(x,y,z))$, where $g\in \mathbb{SE}(3)$, are quadratic and therefore
\begin{equation*}
  Q=\spn \{ 1,x,y,z, x^2, y^2, z^2, xy, xz, yz\} \simeq \mathbb R^{10} \ .  
\end{equation*}
The same implicit representation space arises for any quadratic surface $\bar F$.

Following directly from the definition, if a one--parameter system of surfaces $F_t$ is contained in the homogeneous space $\mathcal H$, then each surface $F_t$ has its implicit equation in $Q$. Using the formalism of the implicit representation mapping \eqref{repM}, it is the zero set of the function $\phi(g_t)\in Q$. 

In Examples \ref{ex1}, \ref{ex2} the equation $\phi(g)$ has the form \eqref{impF} only multiplied by the factor 
\begin{equation*}
    \lambda(t)=\frac{1}{100 \left(1108 t^2+788 t+517\right)^2} \ .  
\end{equation*}
Such a factor does not influence the envelope computation due to Lemma \ref{thm_envelopechar}. The non--uniqueness of the implicit representation suggests a projective construction. More precisely, the implicit representation space could be the projective space $\mathbb P(Q)$. However, exploiting the differential structure of the set of curves would be more difficult in this case. For this reason, we plan to pursue this approach in future research.

We are now ready to formulate and prove the main proposition. First, let us summarize the representations that we defined. Given a one--parameter system $F_t$ of surfaces, we constructed a suitable space $\mathcal H$, which is a homogeneous space of a Lie group $\mathcal G$. We lifted $F_t$ to the group $\mathcal G$ as $g_t$. We defined a linear space of functions $Q$ together with the smooth mapping $\phi$, such that $\phi(g_t)$ is the implicit function of $F_t$. The representation mapping $\phi$ is typically not injective due to a symmetry of $\bar F$. In fact, $\im(\phi)\simeq \mathcal H$ because every surface of $\mathcal H$ has its equation in $Q$.

\begin{thm}\label{thm1}
Let  $F_t$, $t\in I,$ be a one parameter system of implicitly represented surfaces and $\chi_t$ its characteristic curves. Let $\bar F = \left \{[x,y,z] \in \mathbb R^3; \bar{f}(x,y,z)=0 \right \}$ be an elementary surface and suppose $F_t$ is contained in a homogeneous space $\mathcal H = \{g(\bar F): g \in \mathcal G\}$ for a Lie group $\mathcal G$. Let $g_t$ be the lift of $F_t$ into $\mathcal G$. Then for any $t\in I$ there exists an element $\gamma_t\in \mathfrak g$ of the Lie algebra so that 
\begin{equation}\label{bigF}
    \chi_t = g_t(\chi_\gamma),\text{ where }
    \chi_\gamma=\{[x,y,z]\in \mathbb R^3; \bar{f}(x,y,z) = 0 \land d\phi_{\mathbf 1} (\gamma_t)(x,y,z) = 0\} \ .
\end{equation}    
\end{thm}

\begin{prf}
Recall that $d\phi_{\mathbf 1}$ denotes the tangent mapping to $\phi$ taken at the unit element of $\mathcal G$. Because $\phi(\mathbf 1)=\bar{f}$, the domains and co--domain of this linear tangent mapping are as follows:
\begin{equation*}
   d\phi_{\mathbf 1}:T_{\mathbf 1}\mathcal G\to T_{\bar{f}}Q \ .
\end{equation*}
However, by definition, $T_{\mathbf 1}\mathcal G$ is the Lie algebra $\mathfrak g$, and the space $Q$ is linear and therefore identical to its tangent space at any point. We can therefore write 
\begin{equation*}
    d\phi_{\mathbf 1}:\mathfrak g\to Q \ ,
\end{equation*}
and expression \eqref{bigF} is formally consistent. 

For the sake of the proof, let us rename the fixed value of the parameter in $\eqref{bigF}$ from $t$ to $t_0$. This is because the parameter $t$ must remain free for the system description and its differentiation. First, note that $F_{t_0}=g_{t_0}(\bar{F})$, and therefore the characteristic curve $\chi_{t_0}\subset F_{t_0}$ is an image of some curve on $\bar{F}$ under $g_{t_0}$. We will show that it is, in fact, the intersection of $\bar F$ with the zero set of a function of the form $d\phi_{\mathbf 1} (\gamma_{t_0})$. 

To show this, we define the curve $\gamma(t):=g^{-1}_{t_0}g_t$ for $t\in I$. It is simply a transformation of the lifted curve $g_t$ within the group $\mathcal G$, such that $\gamma(t_0)=\mathbf 1$.
We claim that the element $\gamma'(t_0)$ satisfies the formula \eqref{bigF}. Clearly, $\gamma'(t_0)\in \mathfrak g$ because it is the derivative of a curve in $\mathcal G$ at its value corresponding to the unit element. Moreover, 
\begin{equation*}
  f(x,y,z,t)=\phi (g_t)=\phi(g_{t_0}\gamma(t))  
\end{equation*}
is the implicit function of $F_t$. Differentiating this equation with respect to $t$ at the value $t=t_0$ yields 
\begin{equation*}
    \left . \frac{\partial f(x,y,z,t)}{\partial t}\right |_{t=t_0} =d\phi_{\mathbf 1}(\gamma'(t_0))\circ g_{t_0}^{-1} \ ,
\end{equation*}
which concludes the proof, because the right-hand side equation is precisely the transformation by the mapping $g_{t_0}$ of the implicit equation $d\phi_{\mathbf 1}(\gamma'(t_0))$, which represents the derivative surface. \qed
\end{prf}

\medskip
Although the previous proposition may appear complicated due to its great generality, its meaning its meaning is quite easy to explain. All possible characteristic curves are transformations of intersections of one surface, $\bar{f}(x,y,z) = 0$, and surfaces from the linear space, $\im (d\phi_{\mathbf 1})$, see Figure~\ref{fig:ex33}, right. The dimension of this linear space is often relatively low. Indeed, the dimension is equal to the dimension of the Lie group minus the dimension of the surface symmetries. This result leads to elegant computations in Section \ref{apl}.

\section{Applications to envelope parameterization}
\label{apl}
In this section, we apply the previous theory to several systems of surfaces. First, we focus on moving cones, where we provide an explicit rational parameterization of the envelope and we solve the trimming problem. Next, we consider further examples of one--parameter systems generated from different elementary surfaces by various Lie groups. 

\subsection{Rational parameterization of the envelope of moving cones}
Let us continue the analysis of the parametric systems of cones considered in Example \ref{ex1} and Example \ref{ex2}. However, we generalize the setup to consider the Euclidean isometries of an arbitrary cone of revolution.

\begin{exm}\upshape
 Consider the Lie group $\mathcal G=\mathbb{SE}(3)$, fix the parameter  $r\in \mathbb R_+$, and define
\begin{equation*}
    \bar F=\{[x,y,z]: \underbrace{x^2+y^2-r^2z^2}_{\bar f} = 0\} \ ,
\end{equation*}
which is a cone of revolution with an angle $\alpha = \arctan(r)$ between the axis and the generators.
    
The Lie algebra $\mathfrak g= \mathfrak{se}(3)$, see \eqref{Lalg}, is spanned by the elements 
{\footnotesize
\begin{equation}\label{eq_lieAlgRotations}
\gamma_1=\left(
    \begin{array}{cccc}
     0 & 0 & 0 & 0 \\
     0 & 0 & -1 & 0 \\
     0 & 1 & 0 & 0 \\
     0 & 0 & 0 & 0 \\
    \end{array}
    \right),\ 
\gamma_2=\left(
    \begin{array}{cccc}
     0 & 0 & -1 & 0 \\
     0 & 0 & 0 & 0 \\
     1 & 0 & 0 & 0 \\
     0 & 0 & 0 & 0 \\
    \end{array}
    \right),\ 
\gamma_3=\left(
    \begin{array}{cccc}
     0 & -1 & 0 & 0 \\
     1 & 0 & 0 & 0 \\
     0 & 0 & 0 & 0 \\
     0 & 0 & 0 & 0 \\
    \end{array}
    \right), 
\end{equation}
}
which correspond to rotations about the three axes, and
{\footnotesize
\begin{equation}\label{eq_lieAlgTranslations}
\gamma_4=
    \left(
                \begin{array}{cccc}
                 0 & 0 & 0 & 1 \\
                 0 & 0 & 0 & 0 \\
                 0 & 0 & 0 & 0 \\
                 0 & 0 & 0 & 0 \\
                \end{array}
    \right),\  
\gamma_5=\left(
                    \begin{array}{cccc}
                     0 & 0 & 0 & 0 \\
                     0 & 0 & 0 & 1 \\
                     0 & 0 & 0 & 0 \\
                     0 & 0 & 0 & 0 \\
                    \end{array}
                    \right), \ 
\gamma_6=\left(
                        \begin{array}{cccc}
                         0 & 0 & 0 & 0 \\
                         0 & 0 & 0 & 0 \\
                         0 & 0 & 0 & 1 \\
                         0 & 0 & 0 & 0 \\
                        \end{array}
                        \right),
\end{equation}
}
which correspond to translations along the axes. Direct computation yields the implicit representation mapping $\phi$, and its tangent mapping is given by
\begin{equation}\label{comp}
\begin{array}{rcl}
d\phi_{\mathbf 1} (\gamma_1)&=&2 (1+r^2)yz,\\
 d\phi_{\mathbf 1} (\gamma_2)&=&2 (1+r^2)xz,\\
d\phi_{\mathbf 1} (\gamma_3)&=&0,\\
d\phi_{\mathbf 1} (\gamma_4)&=&-2x,\\
 d\phi_{\mathbf 1} (\gamma_5)&=&-2y,\\
d\phi_{\mathbf 1} (\gamma_6)&=&2r^2z.
\end{array}
\end{equation}

Let us explain how these results are computed for the case of $\gamma_1$. It is sufficient to consider a curve in the group $\mathcal G$ representing the tangent element $\gamma_1$. For example, if we take
{\footnotesize
\begin{equation*}
\gamma(t)=\left(
    \begin{array}{cccc}
     1 & 0 & 0 & 0 \\
     0 & \cos(t) & -\sin(t) & 0 \\
     0 & \sin(t) & \cos(t) & 0 \\
     0 & 0 & 0 & 1 \\
    \end{array} 
    \right)\in \mathbb{SE}(3) \ ,
\end{equation*}
}
we clearly have $\gamma'(0)=\gamma_1$. To obtain $d\phi_{\mathbf 1} (\gamma_1)$, we just differentiate $\bar{f}\circ\gamma(t)^{-1}$ at $t=0$. 

The fact that $d\phi_{\mathbf 1} (\gamma_3)=0$ indicates that $\bar F$ is symmetric with respect to rotations about the $z$-axis. The homogeneous space $\mathcal H$ is of dimension $5$ and is isomorphic to the quotient $\mathbb{SE}(3)/\mathbb S^1$.
\end{exm}

The previous computations lead to the following elegant result. 
\begin{thm}\label{thm2}
The characteristic curve of any one--parameter rigid motion of a cone of revolution $\bar F=\{[x,y,z]: x^2+y^2-r^2z^2 = 0\}$ is directly isometric to a curve of the form 
\begin{equation}\label{rigE}
    \chi_h=\{[x,y,z]: x^2+y^2-r^2z^2 = 0 \land  h(x,y,z) = 0\},\text{ where }h(x,y,z)\in \spn \{ x, y, z, xz, yz \} \ .
\end{equation}
Consequently, the characteristic curves are always rational in this case. If the motion is rational, then the envelope surface is also rational.
\end{thm}

\begin{prf}
The term {\em rigid motion} of $\bar F$ merely says that the one--parameter system of surfaces is contained in the homogeneous space $\mathcal H$ with the group $\mathbb{SE}(3)$ and the elementary surface $\bar F$. For this space we computed the images of a basis of $\mathfrak{se}(3)$ under the mapping $d\phi_{\mathbf 1}$, see \eqref{comp}. Since $r$ is a positive constant and the mapping $d\phi_{\mathbf 1}$ is linear we obtain $\im(d\phi_{\mathbf 1})= \spn \{ x, y, z, xz, yz \}$ and the expression \eqref{rigE} follows directly from Proposition \ref{thm1}.

The rationality of $\chi_h$ (and consequently of any $\chi_t$) follows from the fact that the surface $h(x,y,z)=0$ contains the vertex of the cone $\bar F$. Moreover, we were able to perform an explicit parameterization. If $h(x,y,z) = k_1 x + k_2 y + k_3 z + k_4 xz + k_5 yz \in \spn \{ x, y, z, xz, yz \}$, then $\chi_h$ can be parameterized as 
\begin{equation}\label{parCone}
\chi_h(u) =
\begin{pmatrix}
-\frac{(u^2-1) \left(k_1 r u^2-2 k_2 r u-k_1 r+k_3
u^2+k_3\right)}{\left(u^2+1\right) \left(k_4 u^2-2 k_5 u-k_4\right)}\\
\frac{2 u
\left(k_1 r u^2-2 k_2 r u-k_1 r+k_3 u^2+k_3\right)}{\left(u^2+1\right) \left(k_4
u^2-2 k_5 u-k_4\right)}\\
-\frac{k_1 r u^2-2 k_2 r u-k_1 r+k_3 u^2+k_3}{r
\left(k_4 u^2-2 k_5 u-k_4\right)}
\end{pmatrix}.   
\end{equation}

Finally, for rational systems $F_t$, the lifting $g_t$ is rational in $t$ as well as the coefficients $k_i$. Consequently, $g_t(\chi_h(u))$ is a rational parameterization of the envelope surface. 
 \qed
\end{prf}

\medskip
We were able to compute the rational envelope parameterization  $g_t(\chi(u))$ for our running Example \ref{ex1}. Its analytical expression is too long, but we used it to plot the envelope figures through this paper. It is also very useful in analyzing the case of a moving trimmed cone which is crucial for the applications in the CNC manufacturing.  

\begin{figure}[hbt]
    \centering
    \includegraphics[height = 4cm]{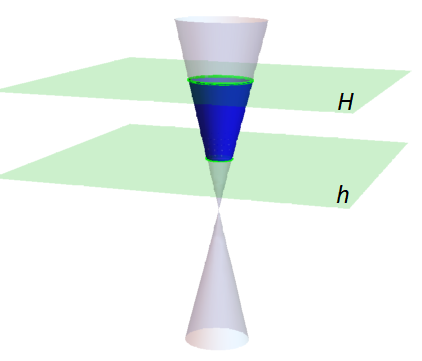}\quad  \includegraphics[height = 3cm]{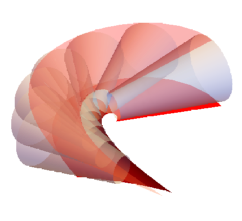}
    \caption{Truncated elementary cone (left) and the exact envelope of its rigid motion (right).}
    \label{fig:truncated}
\end{figure}

\subsection{The trimming problem}
In applications the tool is, of course, bounded. We consider the elementary cone trimmed by two parallel planes perpendicular to the $z$ axis. The resulting envelope then consists of two surface patches (see Fig. \ref{fig:truncated}), each of which is a part of the envelope for the full cone. 

\begin{figure}[tb]
    \centering
 \includegraphics[height = 4.4cm]{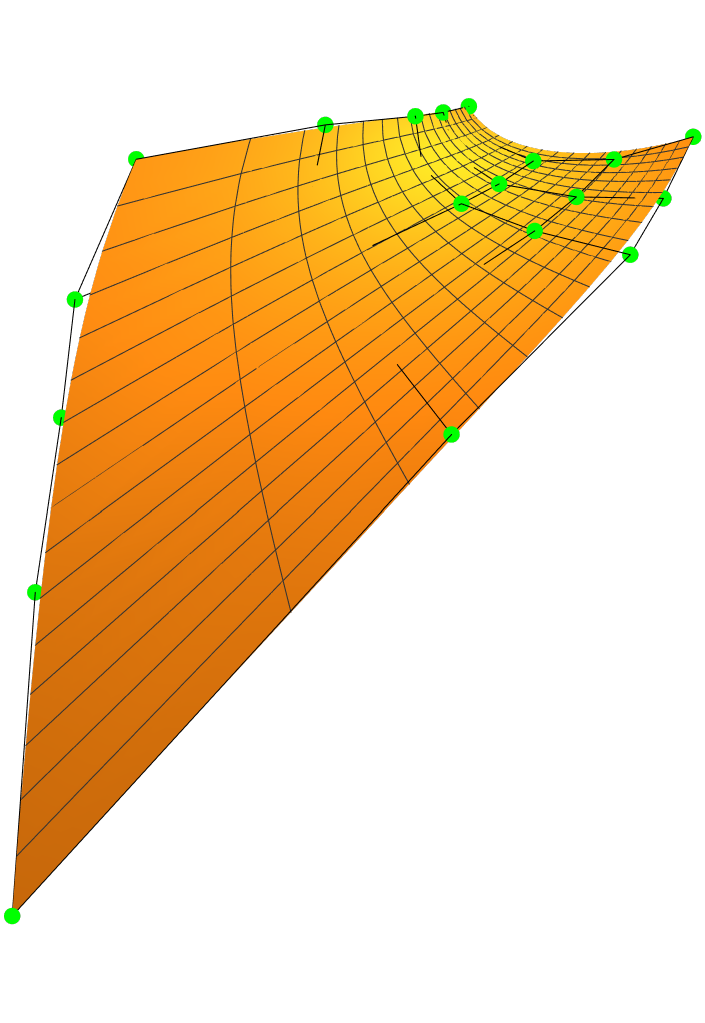}
 \qquad \includegraphics[height = 4cm]{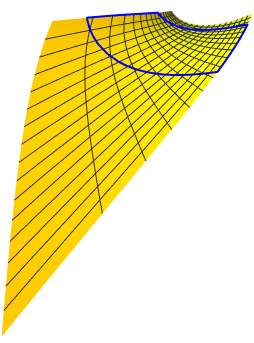}
 \quad \raisebox{.6\height}{\includegraphics[width = 3.5cm]{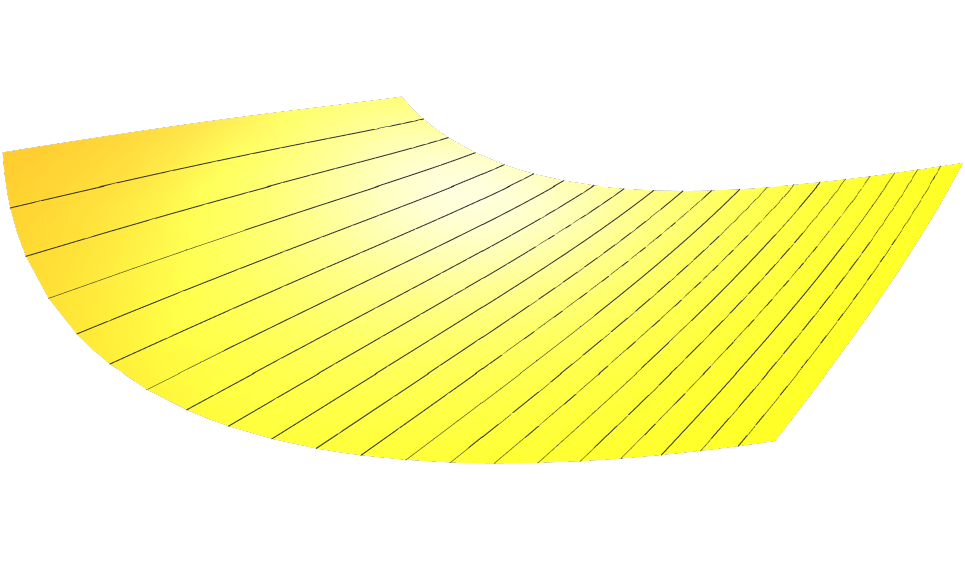}}
    \caption{Surface patch corresponding to the square domain (left) and to the domain trimmed by a B-spline curve (center; right).}
    \label{fig:controlPoints}
\end{figure}

The parameterization \eqref{parCone} can easily be converted into an industrial format. In our example, it can be expressed as a rational Bézier surface, but in general applications some spline motion would be used, resulting in a NURBS envelope surface. Considering some standard domain such as the square $t, u \in [0,1]$, this representation will produce a piece of surface which does not correspond to a trimmed cone, see Fig \ref{fig:controlPoints}, left. 

In order to obtain the correct result, the variable $u$ needs to be bounded depending on the variable $t$. Our explicit rational parameterization allows us to achieve this very efficiently. In the most general case, this leads to the trimmed NURBS format, consisting of a NURBS control points together with NURBS curve trimming, which lives in the domain of the said surface.

\begin{figure}[hbt]
    \centering
    \includegraphics[height = 4cm]{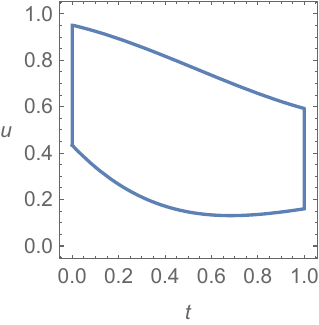}\qquad 
    \includegraphics[height = 4cm]{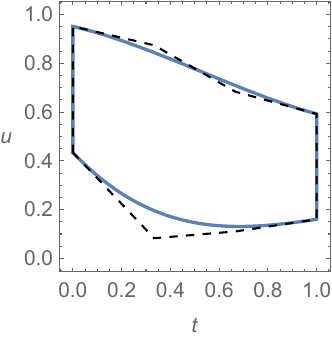}
    \caption{The trimming curve corresponding to one of the envelope patches.}
    \label{fig:trimmingCurve}
\end{figure}

In our example, we obtain the bounding quadratic inequalities for the variable $u$:
\begin{equation*}
    2 < -\frac{k_1 r u^2-2 k_2 r u-k_1 r+k_3 u^2+k_3}{r \left(k_4 u^2-2 k_5 u-k_4\right)} < 5 \ .    
\end{equation*}
They describe the precise domain displayed in Fig. \ref{fig:trimmingCurve}, left. The boundaries can be very efficiently approximated via a functional $C^1$ interpolation by a spline curve displayed in Fig. \ref{fig:trimmingCurve}, right. This way a patch in the {\em trimmed NURBS} format is obtained (see Fig. \ref{fig:controlPoints}, right) which is widely used in the CAD related applications.

\subsection{Examples}
We present three more examples, where the one--parameter systems are obtained by Euclidean and affine transformations of a unit sphere and a paraboloid. These examples are quite general because they represent entire classes of one--parameter systems of surfaces, rather than particular ones.
\begin{exm}
\label{ex3}\upshape
Let us briefly revisit the well--studied case of pipe and canal surfaces, i.e., the envelopes of one--parameter systems of moving spheres with fixed or varying radii, respectively. In both cases, we choose the unit sphere as the elementary surface,
\begin{equation}\label{eq_unitSphere}
    \bar F = \{[x,y,z]: x^2+y^2+z^2-1 = 0\} \ . 
\end{equation}
Any element of the one--parameter system of spheres whose envelope is a pipe surface can be generated from $\bar F$ by the three dimensional group of translations. The corresponding Lie algebra $\mathfrak g$ is spanned by the elements $\gamma_4,\gamma_5$ and $\gamma_6$, see  \eqref{eq_lieAlgTranslations}.

We compute the images of these elements under the tangent mapping $d\phi_1$ and obtain
\begin{align}\label{eq_imagesOfTransl}
    d\phi_1(\gamma_4) &= -2x, \nonumber\\
    d\phi_1(\gamma_5) &= -2y, \\
    d\phi_1(\gamma_6) &= -2z \ . \nonumber
\end{align}
The derivative surface is thus always a plane $h(x,y,z) = k_1x+k_2y+k_3z = 0$ for some $k_1,k_2,k_3 \in \mathbb R$. This plane passes through the center of the sphere and consequently, the characteristic curves of pipe surfaces are the great circles on the spheres, see Fig \ref{fig:ex2}, left.

In order to generate a system of spheres whose envelope is a canal surface, we need to consider one more transformation, which is the expansion mapping of the form 
\begin{equation*}
    \gamma(t)=\left(
        \begin{array}{cccc}
         e^t & 0 & 0 & 0 \\
         0 & e^t & 0 & 0 \\
         0 & 0 & e^t & 0 \\
         0 & 0 & 0 & 1 \\
        \end{array}
        \right) .
\end{equation*}
Computing the image under $d\phi_1$ of the corresponding element in the Lie algebra, we get
\begin{equation*}
    d\phi_{\mathbf 1}(\gamma'(0))=2(x^2+y^2+z^2)=2 \ .
\end{equation*}
Therefore, the derivative surface $h(x,y,z) = k_1x+k_2y+k_3z+k_4 = 0$ for some $k_1,k_2,k_3, k_4 \in \mathbb R$ can be any plane and the characteristic curves for the canal surfaces can be any circles on the spheres, possibly degenerating into a point or an empty set, see Fig \ref{fig:ex2}, right.

\begin{figure}[hbt]
    \centering
   \includegraphics[height=4cm]{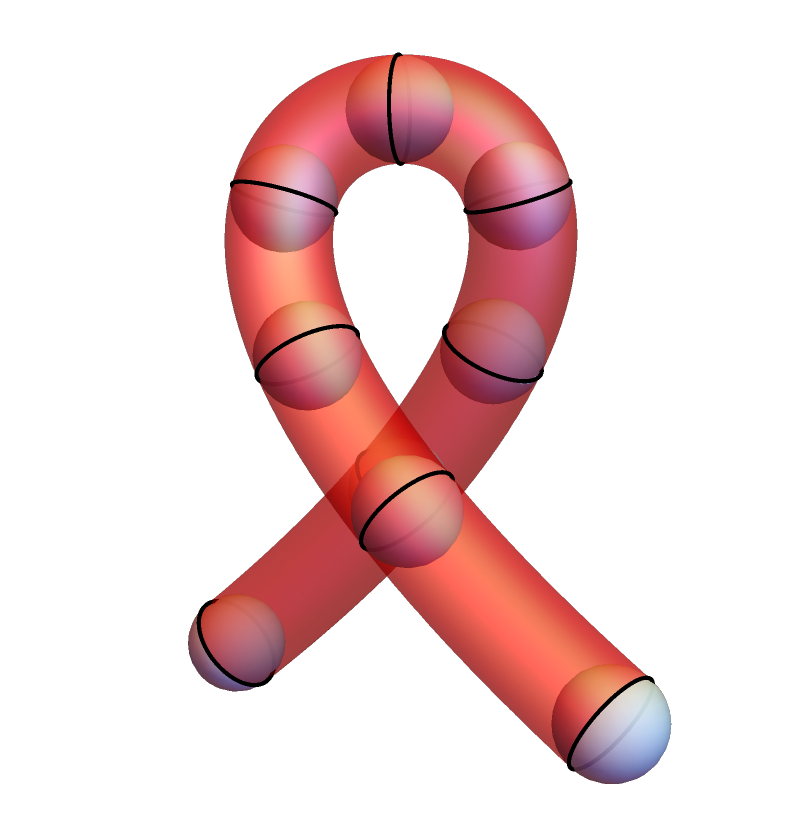}
   \includegraphics[height=4cm]
   {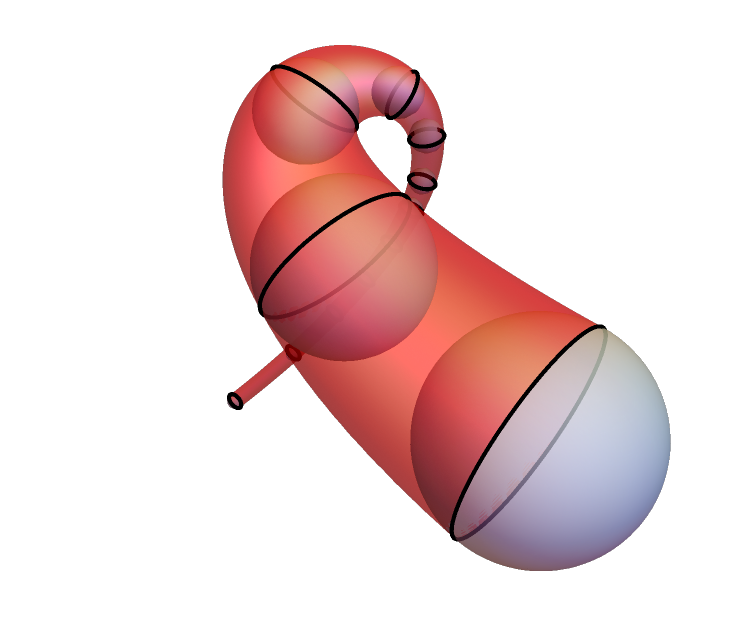}
   \caption{Pipe (left) and canal (right) surfaces as envelopes of one--parameter systems of spheres.}
   \label{fig:ex2}
\end{figure}
\end{exm}

    \begin{exm}
    \label{ex6}\upshape
    Consider again the unit sphere \eqref{eq_unitSphere} as the elementary surface, and a one--parameter system of ellipsoids generated from $\bar F$ by affine transformations. The affine group $\text{Aff}(\mathbb R^3)$ has dimension 12 and can be generated by the three--dimensional group of translations, by expansion mappings of the form
    \begin{equation*}
    \gamma_{e1}(t)=\left(
        \begin{array}{cccc}
         e^t & 0 & 0 & 0 \\
         0 & 1 & 0 & 0 \\
         0 & 0 & 1 & 0 \\
         0 & 0 & 0 & 1 \\
        \end{array}
        \right) , \quad 
        \gamma_{e2}(t)=\left(
        \begin{array}{cccc}
         1 & 0 & 0 & 0 \\
         0 & e^t & 0 & 0 \\
         0 & 0 & 1 & 0 \\
         0 & 0 & 0 & 1 \\
        \end{array}
        \right) , \quad 
        \gamma_{e3}(t)=\left(
        \begin{array}{cccc}
         1 & 0 & 0 & 0 \\
         0 & 1 & 0 & 0 \\
         0 & 0 & e^t & 0 \\
         0 & 0 & 0 & 1 \\
        \end{array}
        \right) ,
\end{equation*}
    and by shear transformations of the form
\begin{align*}
    \gamma_{s1}(t)&=\left(
        \begin{array}{cccc}
         1 & t & 0 & 0 \\
         0 & 1 & 0 & 0 \\
         0 & 0 & 1 & 0 \\
         0 & 0 & 0 & 1 \\
        \end{array}
        \right) , \quad  
    \gamma_{s2}(t)=\left(
        \begin{array}{cccc}
         1 & 0 & t & 0 \\
         0 & 1 & 0 & 0 \\
         0 & 0 & 1 & 0 \\
         0 & 0 & 0 & 1 \\
        \end{array}
        \right) , \quad
    \gamma_{s3}(t)=\left(
        \begin{array}{cccc}
         1 & 0 & 0 & 0 \\
         0 & 1 & t & 0 \\
         0 & 0 & 1 & 0 \\
         0 & 0 & 0 & 1 \\
        \end{array}
        \right) , \\
    \gamma_{s4}(t)&=\left(
        \begin{array}{cccc}
         1 & 0 & 0 & 0 \\
         t & 1 & 0 & 0 \\
         0 & 0 & 1 & 0 \\
         0 & 0 & 0 & 1 \\
        \end{array}
        \right) ,\quad 
    \gamma_{s5}(t)=\left(
        \begin{array}{cccc}
         1 & 0 & 0 & 0 \\
         0 & 1 & 0 & 0 \\
         t & 0 & 1 & 0 \\
         0 & 0 & 0 & 1 \\
        \end{array}
        \right) ,\quad 
    \gamma_{s6}(t)=\left(
        \begin{array}{cccc}
         1 & 0 & 0 & 0 \\
         0 & 1 & 0 & 0 \\
         0 & t & 1 & 0 \\
         0 & 0 & 0 & 1 \\
        \end{array}
        \right) . 
\end{align*}
    By applying $d\phi_{\mathbf 1}$ to the corresponding elements in the Lie algebra $\mathfrak{aff}(\mathbb R^3)$, we obtain
    \begin{equation*}
    \begin{array}{ccc}
        d\phi_{\mathbf 1} (\gamma_{e1}'(0))= -2x^2, & 
        d\phi_{\mathbf 1} (\gamma_{e2}'(0))= -2y^2, &
        d\phi_{\mathbf 1} (\gamma_{e3}'(0))= -2z^2, \\
        d\phi_{\mathbf 1} (\gamma_{s1}'(0))= -2xy, &
        d\phi_{\mathbf 1} (\gamma_{s2}'(0))= -2xz, &
        d\phi_{\mathbf 1} (\gamma_{s3}'(0))= -2yz, \\
        d\phi_{\mathbf 1} (\gamma_{s4}'(0))= -2xy, &
        d\phi_{\mathbf 1} (\gamma_{s5}'(0))= -2xz, &
        d\phi_{\mathbf 1} (\gamma_{s6}'(0))= -2yz. 
    \end{array}
    \end{equation*}
    Including the images of the elements corresponding to translations, see \eqref{eq_imagesOfTransl}, we see that the homogeneous space $\mathcal H$ has dimension 9. All the quadratic functions in $x, y, z$ form a linear space of dimension $10$ but one dimension is lost because of the possible non--zero scalar multiple.  Consequently, any quadratic  surface of the form
    \begin{equation*}
        \{[x,y,z] : h(x,y,z) = 0 \}\ ,\quad \text{where}\quad h(x,y,z) \in \text{Span}\{x,y,z,xy,xz,yz,x^2,y^2,z^2\}\ ,
    \end{equation*}
    can appear as the derivative surface, and the characteristic curves are typically non-rational curves of degree $4$. It is however possible to restrict the allowed transformations and obtain a computationally simpler case. For example Figure \ref{fig:ex3} shows a one--parameter system of ellipsoids obtained from the unit sphere by a subgroup of $\text{Aff}(\mathbb R^3)$ generated by the translations and the shear mapping $\gamma_{s1}$. The derivative surfaces are transformations of a quadratic surface of the form 
    \begin{equation*}
        \{[x,y,z] : k_1x+k_2y+k_3z+k_4xy = 0, k_1,k_2,k_3,k_4 \in \mathbb R \} \ .
    \end{equation*} 

\begin{figure}[hbt]
        \centering
        \includegraphics[height=4cm]{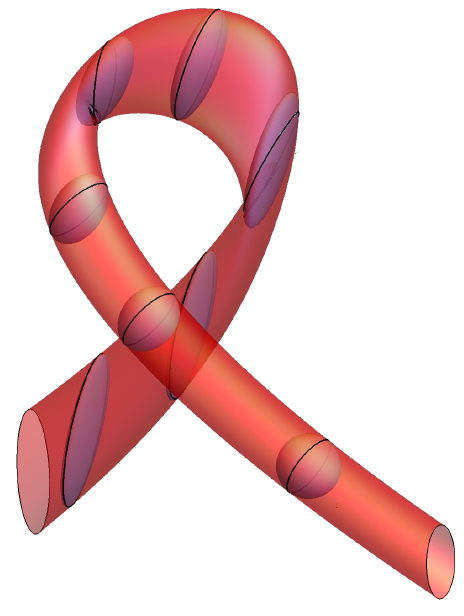}
        \caption{Characteristic curves and the envelope of a one--parameter system of ellipsoids.}
        \label{fig:ex3}
\end{figure}

\end{exm}

\begin{exm}
\upshape
    Finally, consider an elliptical paraboloid
    \begin{equation*}
    \bar F = \{[x,y,z]: \frac{x^2}{a^2}+\frac{y^2}{b^2}-z = 0,\ a,b \neq 0\} \ ,
\end{equation*}
as the elementary surface, and the affine group $\text{Aff}(\mathbb R^3)$. We compute the generators of the homogeneous space $\mathcal H$:

\begin{equation*}
    \begin{array}{lll}
        d\phi_{\mathbf 1} (\gamma_{e1}'(0))= -2ax^2, & 
        d\phi_{\mathbf 1} (\gamma_{e2}'(0))= -2by^2, &
        d\phi_{\mathbf 1} (\gamma_{e3}'(0))= z, \\
        d\phi_{\mathbf 1} (\gamma_{s1}'(0))= -2axy, &
        d\phi_{\mathbf 1} (\gamma_{s2}'(0))= -2axz, &
        d\phi_{\mathbf 1} (\gamma_{s3}'(0))= -2byz, \\
        d\phi_{\mathbf 1} (\gamma_{s4}'(0))= -2bxy, &
        d\phi_{\mathbf 1} (\gamma_{s5}'(0))= x, &
        d\phi_{\mathbf 1} (\gamma_{s6}'(0))= y, \\
        d\phi_{\mathbf 1} (\gamma_4)= -2ax, &
        d\phi_{\mathbf 1} (\gamma_5)= -2by, &
        d\phi_{\mathbf 1} (\gamma_6)= 1,
    \end{array}
    \end{equation*}
    where the last line corresponds to the translations along axes, as in \eqref{eq_lieAlgTranslations}. The homogeneous space has again dimension 9, and the derivative surfaces are affine transformations of a quadratic surface of the form
    \begin{equation*}
        \{[x,y,z] : h(x,y,z) = 0 \}\ ,\quad \text{where}\quad h(x,y,z) \in \text{Span}\{1,x,y,z,xy,xz,yz,x^2,y^2\}\ .
    \end{equation*}

    If we consider the group $\mathbb{SE}(3)$ of Euclidean isometries, the corresponding Lie algebra is generated by $\gamma_1,\gamma_2,\gamma_3$ (corresponding to rotations) and $\gamma_4,\gamma_5,\gamma_6$ (corresponding to translations), see \eqref{eq_lieAlgRotations} and \eqref{eq_lieAlgTranslations}. Computing the images of $\gamma_1,\gamma_2,\gamma_3$ under $d\phi_{\mathbf 1}$ yields
    \begin{equation*}
    \begin{array}{lll}
        d\phi_{\mathbf 1} (\gamma_1)= y+2byz, & 
        d\phi_{\mathbf 1} (\gamma_2)= x+2axz, &
        d\phi_{\mathbf 1} (\gamma_3)= 2(a-b)xy.
    \end{array}
    \end{equation*}
    The homogeneous space has then dimension 6, and the derivative surfaces are isometric to a~quadratic surface of the form
    \begin{equation*}
        \{[x,y,z] : h(x,y,z) = 0 \}\ ,\quad \text{where}\quad h(x,y,z) \in \text{Span}\{1,x,y,xy,xz,yz\}\ .
    \end{equation*}
    We notice that if $a=b$, i.e., if $\bar F$ is a paraboloid of revolution, the dimension drops to 5, because $\bar F$ is stabilized by rotations about the z-axis, see Fig. \ref{ex4}.

\begin{figure}[hbt]
        \centering
        \includegraphics[height=4cm]{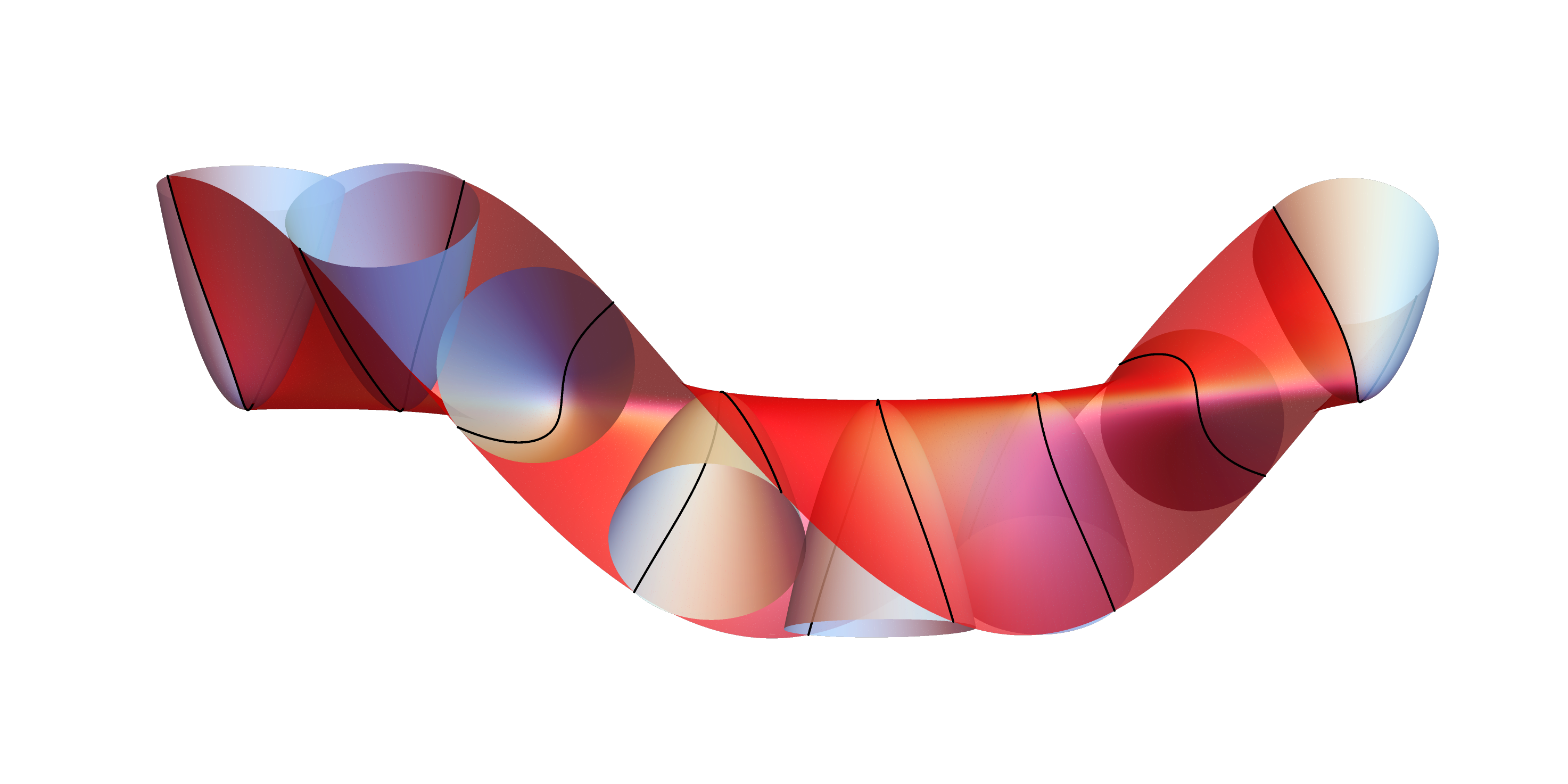}
        \caption{Characteristic curves and the envelope of a one--parameter system of paraboloids of revolution.}
        \label{fig:ex4}
\end{figure}
    
\end{exm}

\section{Conclusion}
We connected the theory of Lie groups and their homogeneous spaces to the description of envelope surfaces. Using the formalism of Lie groups and Lie algebras, we precisely expressed the inherent symmetry and linearity in computing the characteristic curves. For this purpose, we presented several general definitions in Section \ref{Stheory} and formulated the general result in Proposition~\ref{thm1}. 

In Section \ref{apl}, we applied this theory to several examples. In particular, we showed that the characteristic curves and envelopes of cones undergoing rational motion are always rational. While this result is not entirely novel (see the introduction), we provided not only a new, concise proof, but also an explicit method for obtaining the rational parameterization of the envelope.

In the future, we intend to apply the explicit parametrization presented in the proof of Proposition~\ref{thm2} to surface modelling. This will allow us to obtain surfaces that can be manufactured using CNC machining with the flank method and standard conical tools. We also plan to apply the general method to other envelope surfaces. Lastly, we would like to extend the current description to a projective setup.

\bibliography{bibliography}

\end{document}